\documentclass[12pt,twoside,doublespace]{article}

\usepackage{amssymb}
\usepackage{graphicx}

\def\smskip{\par\vskip 5 pt}
\def\QED{\hfill $\Box$\smskip}
\newtheorem{theorem}{Theorem}
\newtheorem{lemma}{Lemma}
\newtheorem{proposition}{Proposition}
\newtheorem{corollary}{Corollary}

\begin{document}

\begin{center}

\vspace{35pt}

{\Large \bf Decomposable Penalty Method }

\vspace{5pt}

{\Large \bf for
Generalized Game Problems

}

\vspace{5pt}

{\Large \bf with Joint Constraints

}

\vspace{35pt}

{\sc I.V.~Konnov\footnote{\normalsize Department of System Analysis
and Information Technologies, Kazan Federal University, ul.
Kremlevskaya, 18, Kazan 420008, Russia.\\ E-mail: konn-igor@ya.ru}}

\end{center}

\begin{abstract}
We consider an extension of a noncooperative game problem where
players have joint binding constraints. In this case, justification of
a generalized equilibrium point needs a reasonable mechanism for
attaining this state. We suggest to combine a  penalty method together
with shares allocation of right-hand sides, which replaces the initial
problem with a sequence of the usual Nash equilibrium problems together
with an upper level variational inequality as a master problem.  We show
  convergence of solutions of these auxiliary penalized problems
to a solution of the initial game problem under weak coercivity conditions.

\end{abstract}


\section{Introduction} \label{s1}

Noncooperative games with joint (binding) constraints date back
to early works 
\cite{Deb52}--\cite{Ich83}  and these problems
are investigated in many recent works; see e.g. \cite{KU00}--\cite{Fuk11}  and
the references therein. Let us consider the
{\em generalized $l$-person noncooperative
game},  where the $i$-th player has its particular strategy set $X_{i} \subseteq
\mathbb{R}^{n_{i}}$  and a payoff (utility) function $f_{i} : X \to
\mathbb{R}$ with
$$
 X=X_{1} \times \dots \times X_{l}, \ n= \sum \limits _{i=1}^{l}
n_{i}.
$$
Besides, all the players together with the above utility functions and strategy sets
have the joint constraint set
$$
  Y  =\left\{x \in \mathbb{R}^{n} \ \vrule \ \sum _{i=1}^{l} h_{i}(x_{i}) \leq b \right\},
$$
where $x= (x_{1},\ldots,x_{l})^{\top}$, $h_{i}(x_{i})= (h_{i1}(x_{i}),\ldots,h_{im}(x_{i}))^{\top}$,
$h_{ij} : \mathbb{R}^{n_{i}} \to \mathbb{R}$, $j=1,\ldots,m$,
$i=1,\ldots,l$ are given functions, $b \in \mathbb{R}^{m}$ is a
fixed vector in $\mathbb{R}^{m}$, as an addition to the set $X$.
That is, they have the common feasible set
\begin{equation} \label{eq:1.5}
 D=X \bigcap Y,
\end{equation}
A point $x^{*}= (x^{*}_{1},\ldots,x^{*}_{l})^{\top}\in D$ is said to
be an equilibrium point for this game, if
\begin{equation} \label{eq:1.6}
    f_{i}(x^{*}_{-i},y_{i})    \leq
     f_{i}(x^{*})   \quad  \forall (x^{*}_{-i},y_{i}) \in D, \ i=1,\ldots,l;
\end{equation}
where we set
$(x_{-i},y_{i})=(x_{1},\ldots,x_{i-1},y_{i},x_{i+1},\ldots,x_{l})$.
The usual noncooperative game corresponds to the case $D=X$, then (\ref{eq:1.6})
is equivalent to the well known Nash equilibrium concept.
Examples of applications of generalized noncooperative
games with joint constraints can be e.g. found in \cite{KU00,CKK04,FK07,PSF08,Kon16c}.
The main question consists in implementation of this
constrained equilibrium concept within the custom noncooperative
framework, where players are independent and make their choices simultaneously.
In fact, the presence of the binding constraints requires certain treaties or
concordant actions of the players, thus contradicting the above assumptions.
These drawbacks were noticed and discussed e.g. in
\cite{Ich83,Har91,Fuk11}.

Rather recently, the right-hand side decomposition approach was suggested for
variational inequalities with binding constraints in \cite{Kon14a}.
Its extension to generalized noncooperative games was given in  \cite{Kon16c}.
Within this approach, the initial problem is treated as a two-level one by using
a share allocation procedure, which leads to a set-valued variational inequality
as a master problem. In  \cite{Kon16c}, a usual noncooperative game problem is solved
at the lower level.
Further,  a decomposable dual regularization (penalty) method that
deals with a single-valued approximation of the master problem for each fixed share allocation
was suggested. Application of this method for production problems with common pollution regulation was described in
\cite{AGK18a}.

In this paper, we suggest to apply a direct decomposable penalty method to generalized noncooperative games,
which involve share allocation variables. We show
 convergence of solutions of auxiliary penalized problems
to a solution of the initial problem under weak coercivity conditions.
Besides, the implementation of this method enables one to solve each
auxiliary penalized problem as a two-level one with the usual Nash equilibrium problem
at the lower level.


\section{Decomposition via Shares Allocation} \label{s2}

We first fix our basic assumptions.

{\bf (A1)} {\em Each strategy set $X_{i}\subseteq
\mathbb{R}^{n_{i}}$ is convex and closed and each
utility function $f_{i}$ is concave in its
$i$-th variable $x_i$ and continuous for $i=1,\ldots,l$. Also,
$h_{ij} : \mathbb{R}^{n_{i}} \to \mathbb{R}$, $j=1,\ldots,m$,
$i=1,\ldots,l$ are convex functions, and the
common feasible set $D $ is nonempty.}

Following the Nikaido-Isoda approach from \cite{NI55}, we consider the normalized
equilibrium problem (EP for short) of
finding a point $x^{*}= (x^{*}_{1},\ldots,x^{*}_{l})^{\top}\in D$ such that
\begin{equation} \label{eq:2.1}
 \Phi(x^{*},y)\geq 0 \quad \forall y \in D,
\end{equation}
where
$$
 \Phi (x,y)=   \Psi (x,x)-  \Psi (x,y), \
 \Psi (x,y) =  \sum _{i=1}^{l} f_{i} (x_{-i},y_{i});
$$
its solutions are called {\em normalized equilibrium points}.
From the above assumptions it follows that $\Phi : X \times X
\rightarrow \mathbb{R}$  is an
equilibrium bi-function, i.e., $\Phi (x,x)=0$ for every $x \in X$,
besides, $\Phi (x, \cdot )$ is convex  for
each $x \in X $ and $\Phi (\cdot, \cdot )$ is continuous.
It should be noted that (\ref{eq:2.1}) implies (\ref{eq:1.6}).
In other words, each normalized
equilibrium point is a generalized Nash equilibrium point, but the
reverse assertion is not true in general. But in case
$D=X$,  (\ref{eq:1.6}) and (\ref{eq:2.1}) become equivalent.
We take a suitable
coercivity condition and obtain the existence result from \cite[Theorem 3.1]{Kon15}.

For a function $\mu : \mathbb{R}^{n} \rightarrow \mathbb{R}$
and a number $r$, we define the level set
$$
    B_{r}= \{ x \in X \ | \ \mu (x) \leq r \}.
$$
We say that the function $\mu : \mathbb{R}^{n} \rightarrow \mathbb{R}$
is weakly coercive with respect to the set $X \subseteq \mathbb{R}^{n}$  if there
exists a number $r$ such that the set $B_{r}$
is nonempty and bounded.

\noindent {\bf (C1)} {\em There exist a lower semicontinuous and
convex function $\mu : X \rightarrow \mathbb{R}$, which is weakly
coercive with respect to the set $D$, and a number $r$,  such that,
for any point $x \in D \setminus B_{r}$ there is a point $z \in D$
with}
$$
\min\{ \Phi(x, z), \mu (z)-\mu (x)\} <0 \quad \mbox{\em and}
\quad\max\{ \Phi(x, z), \mu (z)-\mu (x)\} \leq 0.
$$


\begin{proposition} \label{pro:2.1}
If {\bf (A1)} and {\bf (C1)} are fulfilled, then problem (\ref{eq:2.1}) has a
solution.
\end{proposition}

After proper specialization of the inequalities in {\bf (C1)}
we can somewhat strengthen the above assertion.

\noindent {\bf (C2)} {\em There exist a lower semicontinuous and
convex function $\mu : X \rightarrow \mathbb{R}$, which is weakly
coercive with respect to the set $D$, and a number $r$,  such that,
for any point $x \in D \setminus B_{r}$ there is a point $z \in D$
with}
$$
\mu (z) \leq \mu (x) \quad \mbox{\em and} \quad \Phi(x, z)<0.
$$

\begin{corollary} \label{cor:2.1}
If {\bf (A1)} and {\bf (C2)} are fulfilled, then problem (\ref{eq:2.1}) has a
solution,  and all the solutions are contained in $D \bigcap B_{r}$.
\end{corollary}

These coercivity conditions {\bf (C1)} and {\bf (C2)} clearly hold
if $D$ (or $X$) is bounded. Then we can take $\mu (x)=\|x\|$ and
choose $r$ large enough so that $D \subset B_{r}$.

Let us now introduce the set of  partitions
 of the right-hand side common constraint vector $b$:
$$
\tilde U= \left\{u \in \mathbb{R}^{ml} \ \vrule \ \sum
\limits_{i=1}^{l}u_{i}  = b \right\}.
$$
where $u= (u_{1},\ldots,u_{l})^{\top}$, $u_{i} \in \mathbb{R}^{m}$,
$i=1,\ldots,l$. Here  $u_{i}$ determines the share of the $i$-th player.

Given a partition $u \in \tilde U$, we can consider the
parametric EP:  Find a point $x(u)=(x_{1}(u),\ldots,x_{l}(u))^{\top}
\in D(u)$ such that
\begin{equation} \label{eq:2.2}
 \Phi(x^{*},y)\geq 0 \quad \forall y \in D(u),
\end{equation}
where
$$
 D(u)=D_{1}(u_{1}) \times \dots \times D_{l}(u_{l}),
$$
$D_{i}(u_{i}) = \{ x_{i} \in X_{i}  \ \vrule \ h_{i}(x_{i})
\leq  u_{i} \}$, $i=1,\ldots,l$. Clearly, (\ref{eq:2.2}) is
equivalent to the parametric {\em Nash equilibrium problem} (NEP):
\begin{equation} \label{eq:2.3}
   f_{i}(x_{-i}(u),y_{i}) \leq   f_{i}(x(u))
  \quad  \forall y_{i} \in D_{i}(u_{i}), i=1,\ldots,l.
\end{equation}

If all the optimal shares $u_{i}$, $i=1,\ldots,l$ of players are known,
the constrained problems (\ref{eq:1.6}) and (\ref{eq:2.1}) reduce
to NEPs. Hence, it seems worthwhile to insert
an additional upper control level for finding the optimal shares.
In optimization, this approach is known as the right-hand side
(Kornai-Liptak) decomposition method; see \cite{KL65}.

Following this approach we notice that under certain regularity condition
system (\ref{eq:2.3}) can be replaced with the corresponding system of primal-dual
optimality conditions: Find a pair $(x(u),v(u)) \in X \times
\mathbb{R}_{+}^{ml}$ such that
\begin{eqnarray}
&&  f_{i}(x(u))-f_{i}(x_{-i}(u),y_{i}) \nonumber\\
&&  + \langle v_{i}(u),
h_{i}(y_{i}) -h_{i}(x_{i}(u))\rangle
                \geq 0 \quad \forall y_{i} \in  X_{i}, \label{eq:2.4}\\
&&
 \langle u_{i} - h_{i}(x_{i}(u)), v_{i}-v_{i}(u) \rangle \geq
 0,
  \quad \forall v_{i} \in \mathbb{R}^{m}_{+},  \label{eq:2.5}\\
 && \mbox{for} \ i=1,\ldots,l; \nonumber
\end{eqnarray}
where $v(u)=(v_{1}(u),\ldots,v_{l}(u))^{\top}$.  In this system
the first relations (\ref{eq:2.4}) are rewritten equivalently as
\begin{eqnarray}
&& \displaystyle \Phi (x(u), y) \nonumber\\
&& \displaystyle + \sum _{i=1}^{l} \langle v_{i}(u),
      h_{i}(y_{i}) -h_{i}(x_{i}(u))\rangle
                \geq 0 \quad \forall y \in  X. \label{eq:2.4a}
\end{eqnarray}
Clearly, if
$(x(u),v(u))$ solves (\ref{eq:2.4})--(\ref{eq:2.5}) or (\ref{eq:2.4a}), (\ref{eq:2.5}), then $x(u)$ is
a solution to  (\ref{eq:2.2}) or (\ref{eq:2.3}).

We denote by $T(u)$ the set of all the solution points $-v(u)$,
creating the image of the set-valued mapping $T$. This enables us to
define the variational inequality (VI): Find a point $u^{*} \in
\tilde U$ such that
\begin{equation} \label{eq:2.6}
  \exists t^{*} \in T(u^{*}),
  \ \langle t^{*}, u - u^{*} \rangle   \geq  0, \quad
            \forall  u \in \tilde U.
\end{equation}
 Then it was shown in \cite[Theorem 4.1]{Kon16c} that just the {\em master VI} (\ref{eq:2.6}) yields the
optimal shares of common constraints among players.


\begin{proposition} \label{pro:2.2}
Suppose  {\bf (A1)} is fulfilled. If a  point $u^{*}$ solves VI (\ref{eq:2.6}), the corresponding
point $x(u^{*})$ in (\ref{eq:2.4})--(\ref{eq:2.5}) is a solution
of problem (\ref{eq:1.5}),(\ref{eq:2.1}).
\end{proposition}

We conclude that VI (\ref{eq:2.6}) related to the
parametric problems (\ref{eq:2.2}) or (\ref{eq:2.3})
enables us to find a solution to the initial
generalized noncooperative game. Hence, the two-level
procedure gives a suitable  regulation mechanism for these game problems.
However, this approach has clear drawbacks: $T(u)$ can be empty
for some feasible partitions, besides, $T$ is set-valued in general,
and this fact reduces the number of methods applicable to
solution of  VI (\ref{eq:2.6}). Therefore, this approach
needs certain modifications.


\section{Decomposable Penalty Method} \label{s3}

We start our description of the approach from the simple transformation of the joint constraint
set $Y$ by inserting auxiliary variables:
$$
  Y=\left\{x \in \mathbb{R}^{n} \ \vrule \ \exists u \in \mathbb{R}^{ml},  \sum \limits_{i=1}^{l} u_{i}  =
  b, \ h_{i}(x_{i})  \leq  u_{i}, i=1,\ldots,l \right\},
$$
where $u= (u_{1},\ldots,u_{l})^{\top}$, $u_{i} \in \mathbb{R}^{m}$,
$i=1,\ldots,l$. These variables $u_{i}$ as above determine a partition
 of the right-hand side vector $b$, i.e. give explicit shares of players.
In principle, some additional reasonable restrictions can be imposed on the shares
 $u$, such as $u_{i} \leq b$ or/and $u_{i} \geq \mathbf{0}$ for $i=1,\ldots,l$.
Hence, we define the set of feasible partitions as follows:
$$
U= \left\{u \in U_{0} \ \vrule \ \sum
\limits_{i=1}^{l}u_{i}  = b \right\}.
$$
where $U_{0} \subseteq \mathbb{R}^{ml}$ is a set of these optional additional restrictions
such that for each $x \in D$ there exists $u \in U$,  $h_{i}(x_{i})  \leq  u_{i}$, $i=1,\ldots,l $.

We can now separate the constraints and first consider the auxiliary
penalty problem: Find a pair $w(\tau)= (x(\tau),u(\tau)) \in X
\times U$, $\tau >0$ such that
\begin{equation} \label{eq:3.1}
\Phi_{\tau}(w(\tau), w)=\Phi(x(\tau), x)+\tau [P(w)-P(w(\tau))] \geq
0 \quad \forall w= (x,u) \in X \times U,
\end{equation}
where
$$
P(w)= \sum \limits_{i=1}^{l} P_{i}(w_{i}),
$$
and each $P_{i}(w_{i})=P_{i}(x_{i},u_{i})$ is a general
penalty function for the set
$$
W_{i}=\left\{w_{i}=(x_{i},u_{i}) \in \mathbb{R}^{n_{i}}\times \mathbb{R}^{m}\ \vrule \  h_{i}(x_{i})  \leq  u_{i} \right\}.
$$
We will define these functions as follows:
$$
P_{i}(x_{i},u_{i})=\varphi(h_{i}(x_{i})-u_{i}), \quad i=1,\ldots,l,
$$
where
$\varphi : \mathbb{R}^{m} \to \mathbb{R}_{+}$ is
a convex differentiable and isotone function such that
$$
\varphi(v) \left \{
\begin{array}{cl}
=0, &\quad \mbox{if} \ v \leq \mathbf{0}, \\
>0, & \quad \mbox{otherwise}.
\end{array}  \right.
$$
It follows that
$$
P_{i}(x_{i},u_{i}) \left \{
\begin{array}{cl}
=0, & \quad (x_{i},u_{i}) \in W_{i}, \\ >0, &\quad (x_{i},u_{i}) \notin W_{i};
\end{array}  \right.
$$
for $i=1,\ldots,l$. We recall that the function $ \varphi: \mathbb{R}^{m} \rightarrow
 \mathbb{R} $ is called {\em isotone}, if for
 any points $ u, v$, $u \geq v $ it holds that $ \varphi (u) \geq
 \varphi (v) $. The most popular and simple choice is
$$
\varphi(v)= 0.5 \| [v]_{+}\|^{2},
$$
where $[v]_{+}$ denotes the projection of $v$ onto the non-negative orthant
$\mathbb{R}^{m}_{+}$.
 Then each penalty function $P_{i}$ is convex
and differentiable for $i=1,\ldots,l$.
We observe that decomposable penalty methods were suggested for
 separable convex  optimization problems in \cite{Raz67,Umn75}.

We now show that problem (\ref{eq:3.1}) can be replaced with the following:
Find a pair $w(\tau)= (x(\tau),u(\tau)) \in X \times U$, such that
\begin{eqnarray}
&& \displaystyle \Phi(x(\tau), x)+ \tau \sum \limits_{i=1}^{l} \langle \varphi'(h_{i}(x_{i}(\tau))-u_{i}(\tau)), h_{i}(x_{i})
-h_{i}(x_{i}(\tau))\rangle   \nonumber\\
&& \displaystyle + \tau \sum \limits_{i=1}^{l} \langle \varphi'(h_{i}(x_{i}(\tau))-u_{i}(\tau)), u_{i}(\tau)
-u_{i}\rangle  \geq
0 \quad \forall w= (x,u) \in X \times U.  \label{eq:3.2}
\end{eqnarray}
Let $w(\tau)$ solve (\ref{eq:3.2}). Then (\ref{eq:3.1}) holds due to the
 convexity of the function $\varphi$. In fact, we have
\begin{eqnarray*}
P_{i}(x_{i},u_{i})-P_{i}(x_{i}(\tau),u_{i}(\tau)) \geq &&
\langle \varphi'(h_{i}(x_{i}(\tau))-u_{i}(\tau)), h_{i}(x_{i})-h_{i}(x_{i}(\tau))\rangle \\
&& -\langle \varphi'(h_{i}(x_{i}(\tau))-u_{i}(\tau)), u_{i}-u_{i}(\tau)\rangle,
\end{eqnarray*}
for $i=1,\ldots,l$, and (\ref{eq:3.2}) implies (\ref{eq:3.1}).

 Conversely, let $w(\tau)= (x(\tau),u(\tau))$ be a solution of  problem
 (\ref{eq:3.1}).  Then we can temporarily set $\phi(x)=\Phi(x(\tau), x)$  and obtain
$$
\phi(x)-\phi(x(\tau))+\tau [P(w)-P(w(\tau))] \geq
0 \quad \forall w= (x,u) \in X \times U,
$$
i.e. $w(\tau)$ is a solution of the optimization problem.
Applying now Proposition 4 in \cite{Kon19c} we conclude that
$w(\tau)$ solves (\ref{eq:3.2}).

 In turn, problem (\ref{eq:3.2}) is clearly equivalent to  the system:
Find a pair $w(\tau)= (x(\tau),u(\tau)) \in X \times U$, such that
\begin{eqnarray}
&& \displaystyle \Phi(x(\tau), x)+ \tau \sum \limits_{i=1}^{l} \langle \varphi'(h_{i}(x_{i}(\tau))-u_{i}(\tau)), h_{i}(x_{i})
-h_{i}(x_{i}(\tau))\rangle \geq 0 \quad \forall x \in X,  \label{eq:3.3}\\
&& \displaystyle  \sum \limits_{i=1}^{l} \langle \varphi'(h_{i}(x_{i}(\tau))-u_{i}(\tau)), u_{i}(\tau)-u_{i}\rangle  \geq 0 \quad \forall u \in  U.  \label{eq:3.4}
\end{eqnarray}
We now collect the obtained properties.


 \begin{lemma} \label{lm:3.1} Let the conditions in  {\bf (A1)}  be fulfilled.  Then
 problems (\ref{eq:3.1}), (\ref{eq:3.2}), and (\ref{eq:3.3})--(\ref{eq:3.4}) are equivalent.
\end{lemma}

Given a point $u \in \mathbb{R}^{ml}$, we can solve only problem (\ref{eq:3.3}) in $x$, which is to
find $x(u)\in X$ such that
\begin{equation} \label{eq:3.5}
   \Phi(x(u), x)+ \tau \sum \limits_{i=1}^{l} \langle \varphi'(h_{i}(x_{i}(u))-u_{i}), h_{i}(x_{i})
-h_{i}(x_{i}(u))\rangle \geq 0 \quad \forall x \in X.
\end{equation}
Let $X(u)$ denote the whole solution set of this problem. For each $x(u)\in X(u)$ we set
\begin{equation} \label{eq:3.6}
  g(u)= (g_{1}(u),\ldots,g_{l}(u))^{\top}, \ \mbox{where} \ g_{i}(u)= -\varphi'(h_{i}(x_{i}(u))-u_{i}), \quad i=1,\ldots,l.
\end{equation}
Thus we can define the mapping value
$$
G(u)   =\left\{ g(u) \ | \ x(u)\in X(u) \right\}.
$$
Bearing in mind (\ref{eq:3.4}), we now define the VI: Find a point $u^{*} \in
U$ such that
\begin{equation} \label{eq:3.7}
  \exists g(u^{*}) \in G(u^{*}),
  \ \langle g(u^{*}), u - u^{*} \rangle   \geq  0, \quad
            \forall  u \in U.
\end{equation}


\begin{proposition} \label{pro:3.1}
Suppose  {\bf (A1)} is fulfilled.

(i) If a point $u^{*}$ solves VI (\ref{eq:3.7}), then there exists a point $x^{*}=x(u^{*}) \in X(u^{*})$ such that
$g(u^{*})$ is defined in (\ref{eq:3.6}) at $u= u^{*}$
and that the pair $w^{*}= (x^{*},u^{*})$ is a solution of problem (\ref{eq:3.1}).

(ii) If a pair $w(\tau)= (x(\tau),u(\tau))$ is a solution of problem (\ref{eq:3.1}), then
the point $u^{*}= u(\tau)$ solves VI (\ref{eq:3.7}).
\end{proposition}
The assertions follow directly from the definitions and Lemma \ref{lm:3.1}.

The next step consists in replacing problem (\ref{eq:3.5})
with the following penalized EP: Find $x(u)\in X$ such
that
\begin{equation} \label{eq:3.8}
 \Phi (x(u), y) + \tau  \sum _{i=1}^{l} \left( P_{i}(y_{i},u_{i})-P_{i}(x_{i}(u),u_{i}) \right)
                     \geq 0 \quad  \forall y \in  X.
\end{equation}
The equivalence of (\ref{eq:3.5})  and  (\ref{eq:3.8}) is proved similarly to  Lemma \ref{lm:3.1}.
However, (\ref{eq:3.8}) is clearly equivalent to the NEP: Find
$x(u)\in X$ such that
\begin{equation} \label{eq:3.9}
  \tilde f_{i}(x_{-i}(u),y_{i}) \leq \tilde f_{i}(x(u)),
  \quad  \forall y_{i} \in X_{i}, \ i=1,\ldots,l;
\end{equation}
where the $i$-th player has the penalized utility function
\begin{equation} \label{eq:3.10}
\tilde f_{i}(x)= f_{i}(x)- \tau P_{i}(x_{i},u_{i});
\end{equation}
cf. (\ref{eq:1.6}). Therefore, $X(u)$ is now also the whole solution set of NEP
(\ref{eq:3.9})--(\ref{eq:3.10}) and we have obtained the basic equivalence result.


\begin{theorem} \label{thm:3.1}
Suppose  {\bf (A1)} is fulfilled.

(i) If a point $u^{*}$ solves VI (\ref{eq:3.7}), then there exists a solution $x^{*}=x(u^{*})$
of NEP (\ref{eq:3.9})--(\ref{eq:3.10}) and $g(u^{*})$ is defined in (\ref{eq:3.6}) at $u= u^{*}$,
such that the pair $w^{*}= (x^{*},u^{*})$ is a solution of problem (\ref{eq:3.1}).

(ii) If a pair $w(\tau)= (x(\tau),u(\tau))$ is a solution of problem (\ref{eq:3.1}), then
the point $u^{*}= u(\tau)$ solves VI (\ref{eq:3.7}) and the point $x(\tau)$ is a solution
of NEP (\ref{eq:3.9})--(\ref{eq:3.10}) at $u= u(\tau)$.
\end{theorem}

We conclude that VI (\ref{eq:3.7}) related to the
parametric NEP (\ref{eq:3.9})--(\ref{eq:3.10}) yields a solution for the penalized
game problem. Hence, we have derived another two-level decomposition method
for the initial generalized noncooperative game and have to indicate its preferences over the
method of Section \ref{s2}.
First of all we observe that NEP (\ref{eq:3.9})--(\ref{eq:3.10}) has a
solution under rather mild assumptions. From Proposition \ref{pro:2.1}
we obtain that this is the case if {\bf (A1)}  holds and the set $X$ is bounded.
Also, we can take a suitable coercivity condition in the unbounded case.

\noindent {\bf (C3)} {\em There exists a point $\tilde x \in X$ such that
$$
\Phi (x^{k}, \tilde x) \to -\infty  \quad \mbox{as} \quad \| x^{k}- \tilde
x\| \to \infty;
$$
for any infinite sequence $\{x^{k}\} \subset X$.}


\begin{proposition} \label{pro:3.2}
If {\bf (A1)} and {\bf (C3)} are fulfilled, then NEP (\ref{eq:3.9})--(\ref{eq:3.10}) has a
solution for each $\tau >0$.
\end{proposition}
{\bf Proof.} It suffices to show that problem (\ref{eq:3.8}) has a
solution. Set $\mu (x)=\| x- \tilde x\|$ and
$$
\tilde \Phi (x, y)=\Phi (x, y) + \tau [ P(y,u)-P(x,u)]
$$
for a fixed $u \in U$. Then
$$
\tilde \Phi (x^{k}, \tilde x) \leq \Phi (x^{k}, \tilde x) + \tau P(\tilde x,u) \to -\infty  \quad \mbox{as} \quad \| x^{k}- \tilde
x\| \to \infty.
$$
This means that {\bf (C3)} implies {\bf (C1)} for problem (\ref{eq:3.8})
and the result follows from Proposition \ref{pro:2.1}.
\QED

For this reason, $G(u)$ is non-empty under usual assumptions
even if the set $D(u)$ is empty. The other preference is that $G$
possesses a strengthened monotonicity property.


\begin{proposition} \label{pro:3.3} Suppose {\bf (A1)} is fulfilled,
the bi-function $\Phi$ is monotone on $X \times X$,  i.e.
$$
 \Phi (x',x'') + \Phi (x'',x') \leq 0;
$$
for each pair of points $x',x'' \in X$, and the gradientmap $ \varphi' $
is co-coercive with constant $\gamma$, i.e.
$$
\langle v'-v'', \varphi'(v')-\varphi'(v'') \rangle
\geq \gamma   \|\varphi'(v')-\varphi'(v'') \|^{2}
$$
for all $v',v'' \in \mathbb{R}^{ml}$. Then the mapping $G$ is
co-coercive with constant $\gamma$.
\end{proposition}
{\bf Proof.}
Take arbitrary points $u', u'' \in \mathbb{R}^{ml}$ and set $x'=x(u')$,
$g'=g(u')$, $v'=h(x')-u'$ and $x''=x(u'')$, $g''=g(u'')$, $v''=h(x'')-u''$. It follows from (\ref{eq:3.5})
that
\begin{eqnarray*}
&& \displaystyle \Phi (x',x'')  - \tau\sum _{i=1}^{l} \langle g'_{i},
h_{i}(x''_{i}) -h_{i}(x'_{i})\rangle \geq 0, \\
&& \displaystyle
 \Phi (x'',x') - \tau \sum _{i=1}^{l} \langle g''_{i}, h_{i}(x'_{i})
-h_{i}(x''_{i})\rangle \geq 0;
\end{eqnarray*}
hence
$$
 \sum _{i=1}^{l} \langle g'_{i}- g''_{i},
         h_{i}(x'_{i})-h_{i}(x''_{i}) \rangle \geq  - [\Phi (x',x'')+ \Phi (x'',x')]/\tau \geq 0
$$
since $\Phi $ is monotone. It follows that
\begin{eqnarray*}
&& \displaystyle 0 \leq \sum _{i=1}^{l} \langle g'_{i}- g''_{i},
         [h_{i}(x'_{i})-u'_{i}]-[h_{i}(x''_{i})-u''_{i}] \rangle
         + \sum _{i=1}^{l} \langle g'_{i}- g''_{i}, u'_{i}-u''_{i} \rangle \\
&& \displaystyle
 =  -\sum _{i=1}^{l} \langle \varphi'(v'_{i})-\varphi'(v''_{i}), v'_{i}-v''_{i} \rangle
         + \langle g'- g'', u'-u'' \rangle \\
 && \displaystyle
 \leq -\gamma \sum _{i=1}^{l}    \|\varphi'(v'_{i})-\varphi'(v''_{i}) \|^{2}
 + \langle g'- g'', u'-u'' \rangle \\
&& \displaystyle
 =  -\gamma\|g'- g'' \|^{2} + \langle g'- g'', u'-u'' \rangle ,
\end{eqnarray*}
therefore, $G$ is co-coercive with constant $\gamma$.
\QED

It is well-known that the gradientmap of a convex differentiable function
is co-coercive with constant $\gamma$ if it satisfies the Lipschitz condition with
constant $1/\gamma$; see e.g. \cite[Chapter I, Lemma 6.7]{GT89}.
It follows that the mapping $G$ is then single-valued
even if this is not the case for the mapping $u \mapsto X(u)$. The assertion of
Proposition \ref{pro:3.3} holds true if we replace the monotonicity of
the bi-function $\Phi$ with the more general monotonicity property of the mapping $F$
 defined by taking the sub-differential in $y$ for $\Phi(x,y)$, i.e.,
$$
  F(x)  =  \left. \frac{\partial \Phi (x,y)}{\partial y} \right |_{y=x};
$$
and following the lines of Proposition 7.2 in \cite{Kon16c}. Proposition \ref{pro:3.3}
shows that VI (\ref{eq:3.7}) admits more  efficient solution methods in comparison with VI (\ref{eq:2.6}); see 
e.g. \cite{Kon16c,GT89}.

We observe that the equivalent transformations from  (\ref{eq:3.1}) to  (\ref{eq:3.2}) and
(\ref{eq:3.7})--(\ref{eq:3.8}) are crucial for the decomposable penalty method.
The solution concept based on (\ref{eq:3.7})--(\ref{eq:3.8})  has a
rather simple and natural interpretation. A system regulator chooses first the penalty parameter
$\tau $. Afterwards, he/she determines the right share allocation vector  $u(\tau)$ by
sending some trial vectors $u$ to players and announcing particular deviation penalty functions,
the players then make proper corrections of their utility functions and determine the
corresponding Nash equilibrium point for each trial vector. Then the system regulator changes
the penalty parameter etc. Application of the other known iterative solution methods to
the above generalized game problems was analyzed in \cite{Kon16c}.
They include in particular the usual penalty methods. It appeared that
implementation of these mechanisms within a noncooperative game framework may meet serious
difficulties; see \cite{Kon16c} for more details.


\section{Convergence of the Penalty Method}\label{s4}

In this section, we intend to substantiate the penalty method
with the auxiliary problem (\ref{eq:3.1}) (or (\ref{eq:3.2})).
First we take the following coercivity condition for EP (\ref{eq:3.1}).

\noindent {\bf (C4)} {\em There exist a lower semicontinuous and
convex function $\eta : X \times U \rightarrow \mathbb{R}$, which is weakly
coercive with respect to the set $X \times U$, and a number $r$,  such that,
for any point $w=(x,u) \in (X \times U)$ such that $\eta (w) > r$ there is a point $w'=(z,v) \in (X\times U)$
with}
$$
\eta (w') \leq \eta (w) \quad \mbox{\em and} \quad \Phi_{\tau}(w, w')<0.
$$

We note that {\bf (C4)} is a clear adjustment of condition {\bf (C2)}.
For brevity, we define the level set
$$
    E_{r}= \{ w=(x,u) \in (X \times U) \ | \ \eta (w) \leq r \}.
$$


\begin{lemma} \label{lm:4.1} Let the conditions in  {\bf (A1)} and {\bf (C2)} be fulfilled for some $\tau >0$.  Then
 problem (\ref{eq:3.1}) has a
solution,  and all the solutions are contained in $(X \times U) \bigcap E_{r}$.
\end{lemma}
The assertion follows from Corollary \ref{cor:2.1}.

However,  condition {\bf (C4)} is not suitable for verification (cf. {\bf (C3)}) and we will deduce
it from other conditions of form {\bf (C2)}. In general, we follow the approach from \cite{Kon15}.
First we consider the case where the set $U$ is bounded.

\noindent {\bf (C5)}  {\em There exist a lower semicontinuous and
convex function $\mu : X \rightarrow \mathbb{R}$, which is weakly
coercive with respect to the set $X$, and a number $r$,  such that,
for any point $x \in X \setminus B_{r}$ there is a point $z \in D$
with}
$$
\mu (z) \leq \mu (x) \quad \mbox{\em and} \quad \Phi(x, z)<0.
$$


\begin{theorem} \label{thm:4.1}
Suppose that {\bf (A1)} and {\bf (C5)} are fulfilled, the set $U$ is bounded, the sequence
$\{\tau_{k}\}$  satisfies
\begin{equation} \label{eq:4.1}
\{\tau_{k}\} \nearrow +\infty.
\end{equation}
Then:

(i) EP (\ref{eq:2.1}) has a solution;

(ii) EP (\ref{eq:3.1}) has a solution for each  $\tau > 0$
 and all these solutions belong to $B_{r}\times U $;

(iii) Each sequence $\{w(\tau_{k})\}$ of solutions of EP (\ref{eq:3.1})
has limit points, all these limit points belong to $(B_{r}\bigcap D)\times U $
 and are solutions of EP (\ref{eq:2.1}).
 \end{theorem}
{\bf Proof.} We first show that, for any $\tau > 0$, {\bf (C4)} is true
with $\eta (w)=\mu(x)$. Take any $w=(x,u) \in (X \times U) \setminus E_{r}$,
then by {\bf (C5)} there is $z \in D$, $\mu (z) \leq \mu (x)$ such that $\Phi(x, z)<0$.
Since $z \in D$, there exists $v \in U$ such that $P(w')=P(z, v)=0$ and we have
$$
  \Phi _{\tau}(w, w')  =  \Phi (x,z)+\tau [ P( w')-P ( w)] \leq \Phi (x,z)-\tau P ( w)
  < 0.
$$
 Hence, assertion (ii) follows from Lemma \ref{lm:4.1}.

For brevity, we set $w^{k}=w(\tau_{k})=(x^{k},u^{k})$, where $x^{k}=x(\tau_{k})$ and $u^{k}=u(\tau_{k})$.
By (ii), the sequence $\{w^{k}\}$ exists and is bounded.
Therefore, it has limit points. Since $B_{r}$ is
convex and closed, all these limit points  must belong to
$B_{r}\times U $. Let $\bar w=(\bar x,\bar u)$ be an arbitrary limit point of
$\{w^{k}\}$, i.e. $ \{w^{k_{s}} \} \to \bar w$. Then, by
definition,
$$
 0 \leq P(w^{k_{s}})\leq\tau_{k_{s}}^{-1}\Phi(x^{k_{s}},x)+ P(x,u), \quad \forall (x,u) \in X \times U.
$$
Taking $(x,u) \in D \times U $ such that $P(x,u)=0$ and using (\ref{eq:4.1}), we obtain
$$
  0 \leq P(\bar w) \leq \liminf _{s \rightarrow \infty}  P(w^{k_{s}})\leq \limsup _{s \rightarrow \infty}
 \left[\tau_{k_{s}}^{-1}\Phi(x^{k_{s}},x) \right] \leq 0,
$$
i.e. $P(\bar w)=0$ and $\bar x \in D$.
Therefore, for each $x \in D$ we can take $u \in U$ such that $P(w)=P(x, u)=0$ and
obtain
$$
\Phi (x^{k_{s}}, x) -  \tau_{k_{s}} P(w^{k_{s}})
   =\Phi(x^{k_{s}},x)+ \tau_{k_{s}}\left[P(w)-  P(w^{k_{s}})\right]  \geq 0.
$$
It now follows that
$$
 \Phi (\bar x, x) \geq  \limsup _{s \rightarrow +\infty } \Phi (x^{k_{s}},x)
  \geq  \limsup_{s \rightarrow +\infty } \left[\tau_{k_{s}}P(x^{k_{s}}) \right]
  \geq 0.
$$
Therefore $\bar x$ solves EP (\ref{eq:2.1}) and assertion
(iii) is true. Since $\bar x$ exists, assertion (i) is also true. The
proof is complete. \QED

We now give  similar properties in the unbounded case.


\begin{theorem} \label{thm:4.2}
Suppose that {\bf (A1)} and {\bf (C5)} are fulfilled,  the sequence
$\{\tau_{k}\}$  satisfies condition (\ref{eq:4.1}).
Then:

(i) EP (\ref{eq:2.1}) has a solution;

(ii) There exists $\tau' > 0$ such that EP (\ref{eq:3.1}) has a solution for each  $\tau > \tau'$
 and all these solutions belong to $E_{r}\bigcap ( X\times U) $ where $\eta (w)=\max\{\mu(x),P(w)\}$;

(iii) Each sequence $\{w(\tau_{k})\}$ of solutions of EP (\ref{eq:3.1})
has limit points, all these limit points belong to $(B_{r}\bigcap D)\times U $
 and are solutions of EP (\ref{eq:2.1}).
 \end{theorem}
{\bf Proof.} We first show that the function $\eta (w)=\max\{\mu(x),P ( w)\}$ is weakly
coercive with respect to the set $X \times U$. Fix a number $r$ such that the set
$E_{r}\bigcap ( X\times U) $ is nonempty. If it is unbounded,
there exists a sequence $\{w^{k}\}$ such that $w^{k} =(x^{k},u^{k}) \in E_{r}\bigcap ( X\times U)$ such that
$\| w^{k}\| \to \infty $ as $k \to \infty$. Since
the function $\mu(x)$ is weakly coercive with respect to the set $X $, we have
$\| x^{k}\| \leq C \leq \infty$. It follows that
$\| u^{k}\| \to \infty$. Since $u \in  U$, there exists at
least one pair of indices $j$ and $t$ such that $ u^{k_{s}}_{jt} \to
-\infty$ for the corresponding subsequence $\{u^{k_{s}}\}$,
hence $ P_{j}(x^{k_{s}}_{j},u^{k_{s}}_{j}) \to +\infty$, which is
a contradiction.

We now show that there exists $\tau' > 0$ such that {\bf (C4)} is true
for any $\tau > \tau'$. Take any $w=(x,u) \in (X \times U) \setminus E_{r}$, then
$\eta(w)> r$. If $\mu (x)> r$,
then by {\bf (C5)} there is $z \in D$, $\mu (z) \leq \mu (x)$ such that $\Phi(x, z)<0$.
Since $z \in D$, there exists $v \in U$ such that $P(w')=P(z, v)=0$ and we have
$$
  \Phi _{\tau}(w, w')  =  \Phi (x,z)+\tau [ P( w')-P ( w)] \leq \Phi (x,z)-\tau P ( w)
  < 0.
$$
Hence, {\bf (C4)} holds. If $\mu (x) \leq r$, then $P ( w)> r$. Fix a point $z \in D$, then
there exists $v \in U$ such that $P(w')=P(z, v)=0$. Since the set $B_{r}$ is bounded, we have
$$
  \max\limits_{x \in B_{r}} \Phi(x, z) =  d < \infty.
$$
Take any $\tau'=d/r $, then
$$
  \Phi _{\tau}(w, w')  =  \Phi (x,z)+\tau [ P( w')-P ( w)] \leq d-\tau r < 0
$$
if $\tau > \tau'$.
Hence, {\bf (C4)} also holds and assertion (ii) follows from Lemma \ref{lm:4.1}.
Assertions (iii) and (i) are proved as in Theorem \ref{thm:4.1}. \QED


\end{document}